\documentclass[11pt,a4paper]{article}
\usepackage[english]{babel}
\usepackage{amsfonts}
\usepackage{amsthm}
\newtheorem{teo}{Theorem}[section]
\newtheorem{prop}[teo]{Proposition}
\newtheorem{lema}[teo]{Lemma}
\newtheorem{coro}[teo]{Corollary}
\newtheorem{defi}[teo]{Definition}
\theoremstyle{definition}
\newtheorem{rem}[teo]{Remark}
\newtheorem{ejem}[teo]{Example}

\begin{document}

\title{\vspace*{0cm}Unitary Orbits in a Full Matrix Algebra}
\date{}
\author{Gabriel Larotonda}

\maketitle

\abstract{\footnotesize{\noindent The Hilbert manifold $\Sigma$ consisting of positive invertible (unitized) Hilbert-Schmidt operators has a rich  structure and geometry.  The geometry of unitary orbits $\Omega\subset \Sigma$ is studied from the topological and metric viewpoints: we seek for conditions that ensure the existence of a smooth local structure for the set $\Omega$, and we study the convexity of this set for the geodesic structures that arise when we give $\Sigma$  two Riemannian metrics.
}}

\section{Introduction}

In this paper we study the geometry of unitary orbits $\Omega$ in a (Riemannian, infinite dimensional) manifold $\Sigma_{\infty}$, a manifold which is modeled on the full-matrix algebra of Hilbert-Schmidt operators. We investigate necessary and sufficient conditions for these orbits to be analytic submanifolds. We are also concerned with the explicit form of the geodesics in such submanifolds, and whether this submanifolds are convex when embedded in the full space $\Sigma_{\infty}$, or even in the tangent (Euclidean) space (where $\Sigma_{\infty}$ is open). The last results of this paper give a satisfactory characterization of the exponential map of the submanifold $\Omega$ when this set is the unitary orbit of a projection.

The main framework of this paper is the von Neumann algebra ${\sf B}(H)$ of bounded operators acting on a complex, separable Hilbert space $H$.

Throughout,  ${{\sf HS}}$ stands for the bilateral ideal of Hilbert-Schmidt operators  of ${\sf B}(H)$. This ideal is known as a full-matrix algebra \cite{rickart} since any Hilbert-Schmidt operator can be identified with an infinite matrix such that any row (and any column) is square-summable. Recall  \cite{simon} that ${{\sf HS}}$ is a Banach algebra (without unit) when given the norm $\|a\|_{2}=2\;tr(a^*a)^{\frac{1}{2}}$. Inside ${\sf B}(H)$ we consider a certain kind of Fredholm operators, namely
$${\sf H}_{\mathbb  C}=\{ a+\lambda: \;\; a\in {{\sf HS}},\; \lambda\in\mathbb C \},$$
the complex linear subalgebra consisting of Hilbert-Schmidt perturbations of scalar multiples of the identity. Note that this is a complex Hilbert space with the inner product
$$
\left<\alpha+a,\beta+b\right>_{_2}=\alpha\overline{\beta}+4tr(b^*a)
$$
The model space that we are interested in is the real part of ${\sf H}_{\mathbb C}$:
$${\sf H}_{\mathbb R}=\{ a+\lambda: \;a^*=a,\; a\in {{\sf HS}},\; \lambda\in\mathbb R \},$$
which inherits the structure of real Banach space, and with the same inner product, becomes a real Hilbert space. 

\begin{rem}\label{darvuelta}For this inner product, we have (by cyclicity of the trace)
$$
\left<XY,Y^*X^*\right>_{_2}=\left<YX,X^*Y^*\right>_{_2}\quad\mbox{ for any }X,Y\in {\sf H}_{\mathbb C}\mbox{, and  also}
$$
$$
\left<ZX,YZ\right>_{_2}=\left<XZ,ZY\right>_{_2}\quad\mbox{ for }X,Y\in {\sf H}_{\mathbb C}\mbox{ and }Z\in {\sf H}_{\mathbb R}
$$
\end{rem}

\smallskip

We will use ${{\sf HS}}^h$ to denote the closed subspace of self-adjoint Hilbert-Schmidt operators. Inside ${\sf H}_{\mathbb R}$, consider the subset 
$$\Sigma_{\infty}:=\{A>0, A\in {\sf H}_{\mathbb R}\}$$
This is the set of invertible operators $a+\lambda$ such that $\sigma(a+\lambda) \subset (0,+\infty)$, with $a$ self-adjoint and Hilbert-Schmidt, $\lambda \in \mathbb R$. Note that, since $a$ is compact, then $0\in \sigma(a)$, which forces $\lambda>0$. It is apparent that $\Sigma_{\infty}$ is an open set of ${\sf H}_{\mathbb R}$, therefore a real analytic submanifold. For any $p\in\Sigma_{\infty}$, we may thus identify $T_p \Sigma_{\infty}$ with  $ {\sf H}_{\mathbb R}$, and endow this manifold with a (real) Riemannian metric by means of the formula
$$
\left<X,Y\right>_{_p}=\left<p^{-1}X,Yp^{-1}\right>_{_2}=\left<Xp^{-1},p^{-1}Y\right>_{_2}
$$
With this metric $\Sigma_{\infty}$ has nonpositive sectional curvature \cite{av1}; moreover, the curvature tensor is given by the following commutant:
\begin{equation}\label{tensor}
{\sf R}_p(X,Y)Z=-\frac14 \; p\left[  \left[p^{-1}X,p^{-1}Y \right]  ,  p^{-1}Z       \right]
\end{equation}
Covariant derivative is given by the expression
\begin{equation}\label{covariant}
\nabla_X Y=X(Y)-\frac{1}{2}\left( Xp^{-1}Y+Yp^{-1}X  \right)
\end{equation}
where $X(Y)$ denotes derivation of the vector field $Y$ in the direction of $X$ (performed in the ambient space ${\sf H}_{\mathbb R}$). Euler's equation $\nabla_{\dot \gamma}\dot\gamma=0$ reads
\begin{equation}\label{geo}
\ddot\gamma\ - \dot\gamma\gamma^{^{-1}}\dot\gamma=0,
\end{equation}
and the unique geodesic joining $\gamma_{pq}(0)=p$ with $\gamma_{pq}(1)=q\,$ is given by the expression
\begin{equation}\label{curvas}
\gamma_{pq}(t)=p^{\frac{1}{2}}\left( p^{-\frac{1}{2}}q p^{-\frac{1}{2}}\right)^t p^{\frac{1}{2}}
\end{equation}
These curves look formally equal to the geodesics between positive definite matrices (regarded as a symmetric space); this geodesic is unique and realizes the distance: the manifold $\Sigma_{\infty}$ turns out to be complete with this distance. 

\bigskip

\begin{rem}\label{exponencial}
Throughout,  $\|X\|_{_p}^2:=\left<X,X\right>_{_p}$, namely
$$\|X\|_{_p}^2=\|p^{-1/2}Xp^{-1/2}\|_{_2}=\left<X p^{-1},p^{-1}X\right>_{_2}=\left<p^{-1}X,Xp^{-1}\right>_{_2},$$
which is the norm of tangent vectors $X\in T_p\Sigma_{\infty}$. We will use ${\rm exp}_p$ to denote the exponential map of $\;\Sigma_{\infty}$. 

Note that ${\rm exp}_p(V)=p^{\frac12}\;{\rm e}^{\;p^{-\frac12}\,V\,p^{-\frac12}}p^{\frac12},$ but rearranging the exponential series we get the alternate expressions\index{exponential map}\index{_exponentia@${\rm exp_p}$ exponential map at $p$}
$$
{\rm exp}_p(V)=p\;{\rm e}^{p^{-1}V}={\rm e}^{\,Vp^{-1}}p
$$
A straightforward computation also shows that for $p,q\in\Sigma_{\infty}$ we have 
$${\rm exp}_p^{-1}(q)=p^{\frac12}\ln(p^{-\frac12}\,q\,p^{-\frac12})p^{\frac12}$$
\end{rem}

\begin{lema}\label{invar}\index{group!of isometries}
The metric in $\Sigma_{\infty}$ is invariant for the action of the group of invertible elements: if $g$ is an invertible operator in ${\sf H}_{\mathbb C}$, then $I_g(p)=gpg^*$ is an isometry of $\Sigma_{\infty}$.\index{isometry}
\end{lema}
\begin{proof}It follows from Remark \ref{darvuelta}.\end{proof}

\smallskip

\begin{rem}
$\Sigma_{\infty}$ is complete as a metric space due to the following fact, which is also strongly connected with the fact that sectional curvature is nonpositive, see \cite{av1} and \cite{gnegk}:
$$
\|X-Y\|_{_2}\le \mbox{\rm dist}({\rm e}^X,{\rm e}^Y)=\|\ln({\rm e}^{-X/2}{\rm e}^Y{\rm e}^{-X/2})\|_{_2}
$$
This inequality was first shown for the operator (spectral) norm in the paper \cite{cpr1}, and in that context is related with I. Segal's inequality $\|e^{x+y}\|_{\infty}\le \|e^{x/2}e^y e^{x/2}\|_{\infty}$ (see \cite{cpr92} for further details).

The manifold $\Sigma_{\infty}$ is also complete in the following sense: ${\rm exp}_p$ is a diffeomorphism onto $\Sigma_{\infty}$ for any $p$. The reader should be careful with other notions of completeness, because, as C.J. Atkin shows in \cite{atkin1} and \cite{atkin2}, Hopf-Rinow theorem is not valid in the infinite dimensional context.
\end{rem}

\bigskip

\section{Unitary orbits}\label{orbita}

The total manifold can be decomposed as a disjoint union of geodesically convex submanifolds
$$
\Sigma_{\lambda}=\{a+\lambda\in \Sigma_{\infty},\; a\in {{\sf HS}}^h \textnormal{ and }\lambda>0 \textnormal{ fixed }\}
$$
There is a distinguished leaf in the foliation, namely $\Sigma_{1}$, which contains the identity. Moreover, $\Sigma_1=\mbox{\rm exp}({{\sf HS}}^h)$. We will focus on this submanifold since the nontrivial part of the geometry of $\Sigma_{\infty}$ is contained in the leaves \cite{gnegk}. We won't have to deal with the scalar part of tangent vectors, and some computations will be less involved.

\medskip

\subsection{\large The action of the unitary groups ${\sf U}{{\sf H}_{\mathbb C}}$, ${\sf U}{{\sf B}(H)}$}

We are interested in the orbit of an element $1+a\in \Sigma_1$ by means of the action of some group of unitaries. We first consider the group of unitaries of the complex Banach algebra of 'unitized' Hilbert-Schmidt operators. To be precise, let's call
$$
{\sf U}{\sf H}_{\mathbb C}=\{g=1+a:a\in {{\sf HS}},  g^*=g^{-1}\}
$$
The Lie algebra of this Lie group consists of the operators of the form $ix$ where $x$ is a Hilbert-Schmidt, self-adjoint operator
$$
T_1\left({\sf U}{\sf H}_{\mathbb C} \right)=i{\sf HS}^{h}={\sf HS}^{ah}
$$

\begin{rem}
The problem of determining whether a set in $\Sigma_1$ can be given the structure of submanifold (or not) can be translated into the tangent space by taking logarithms; to be precise, note that 
$$
\mbox{\rm exp}(UaU^*)=U\mbox{e}^a U^*
$$
for any $a\in {{\sf HS}}^h$ and any unitary $U$, and that this map is an analytic isomorphism between $\Sigma_1$ and its tangent space. We will state the problem in this context.
\end{rem}

\smallskip

We fix an element $a$ in the tangent space (that is, $a\in {{\sf HS}}^h$) and make the unitary group act {\it via} the map
$$
\pi_a:{\sf U} {\sf H}_{\mathbb C} \to {{\sf HS}}^h\qquad  g\mapsto gag^*
$$

\begin{defi}
Let $S_a$ be the orbit of the element $a\in {\sf HS}^h$ for the action of the Hilbert-Schmidt unitary group, that is $S_a=\pi_a \left( {\sf U} {\sf H}_{\mathbb C} \right)$.
\end{defi}

This raises the question: when is the orbit of a self-adjoint Hilbert-Schmidt operator a submanifold of ${{\sf HS}}^h$? The answer to this question can be partially answered in terms of the spectrum of the fixed operator:

\begin{teo}\label{hs} If the algebra $C^*(a)$ generated by $a$ and $1$ is finite dimensional, then the orbit $S_a\subset {{\sf HS}}^h$ can be given an analytic submanifold structure.
\end{teo}
\begin{proof} A local section for the map $\pi_a$ is a pair $(U_a,\varphi_a)$ where $U_a$ is an open neighborhood of $a$ in ${{\sf HS}}^h$ and $\varphi_a$ is an analytic map from $U_a$ to ${\sf U} {\sf H}_{\mathbb C}$ such that:
\begin{itemize}
\item $\varphi_a(a)=1$
\item $\varphi_a$ restricted to $U_a\cap S_a$ is a section for $\pi_a$, that is 
$$\pi_a\circ\varphi_a\left|_{U_a\cap S_a}\right.=id_{U_a\cap S_a}$$
\end{itemize} 

A section for $\pi_a$ provides us with sufficient conditions to give the orbit the structure of immerse submanifold of ${{\sf HS}}^h$ (see Proposition 2.1 of \cite{as1}). The section $\varphi_a$ can be constructed by means of the finite rank projections in the matrix algebra where $C^*(a)$ is represented. The finite dimension of the algebra is key to the continuity (and furthermore analyticity) of all the maps involved. To fix some notation, set $n=\dim\;C^*(a)$ and $\tau$ an $^*$-isomorphism
$$
\tau: C^*(a)\to \mathbb C\oplus \mathbb C\oplus\cdots\oplus \mathbb C
$$
Consider the set of systems of one-dimensional projections (here $p_i^2=p_i=p_i^*, p_ip_j=0$ for any $i\ne j$):
$$
P_n=\{(p_1,\cdots,p_n)\in {\sf H}_{\mathbb C}^n:\;\sum_{i=1}^n p_i=1\}
$$
Denote $e_{jk}^i\in M_{n_i}(\mathbb C)$ the elementary matrix with $1$ in the $(j,k)$-entry and zero elsewhere, but embedded in the direct sum; take  $p_{jk}^i(X)$ the polynomial which makes $e_{jk}^i=p_{jk}^i(\tau(a))$, and consider the following element in ${{\sf HS}}^h$: ${\rm E}^i_{jk}=p_{jk}^i(a)$

There is a neighborhood $U_a$ of $a$ in ${{\sf HS}}^h$ such that $1-\left[{\rm e}^i_{11}-p_{11}^i(x)\right]^2$ has strictly positive spectrum, because $r(x)=\|x\|\le \|x\|_{_2}$ and ${\sf H}_{\mathbb C}$ is a Banach algebra (here $r(x)$ denotes spectral radius). A straightforward computation shows that the map
$$
\varphi_a(x)=\sum_{i=1}^p \sum_{j=i}^{n_i} p_{j1}^i(x) E_{11}^i\left[ 1-\left( E_{11}^i-p_{11}^i(x)\right)^2\right]^{-\frac{1}{2}}\;E_{1j}^i
$$
is a cross section for $\pi_a$, and it is analytic from $U_a\subset {{\sf HS}}^h\to {\sf U} {\sf H}_{\mathbb C}$ since the $p^i_{jk}$ are multilinear and all the operations are taken inside the Banach algebra ${\sf H}_{\mathbb C}$. \end{proof}

\smallskip

\begin{rem}\label{nose} At first sight, it is not obvious if this strong restriction (on the spectrum of $a$) is necessary for $S_a$ to be a submanifold of ${{\sf HS}}^h$. The main difference with the work done so far by Deckard and Fialkow in \cite{deckard}, Raeburn  in \cite{raeburn}, and Andruchow \emph{et al}. in \cite{as1}, \cite{as2} is that the Hilbert-Schmidt operators (with any norm equivalent to the $\|\cdot\|_{_2}$-norm) are not a $C^*$-algebra. A remarkable byproduct of Voiculescu's theorem \cite{voicu} says that, for the unitary orbit of an operator $a$ with the action of the full group of unitaries of ${\sf B}(H)$, it is indeed necessary that $a$ has finite spectrum. For the time being, we don't know if this is true for the algebra ${\sf B}={\sf H}_{\mathbb C}$. \end{rem}

\smallskip
\index{group!action}

Let's examine what happens when we act with the full unitary group ${\sf U} {\sf B}(H)$ by means of the same action. For convenience let's fix the notation

$${\mathfrak S}_a=\{UaU^*\;:\;U\in {\sf U} {\sf B}(H)\}$$

We will develop an example that shows that the two orbits ($S_a$ and $\mathfrak S_a$) are, in general, not equal when the spectrum of $a$ is infinite.

\begin{ejem}\label{noda} Take $H=l_2(\mathbb Z)$, $S\in {\sf B}(H)$ the bilateral shift ($Se_k=e_{k+1}$). Then $S$ is a unitary operator with $S^*e_k=e_{k-1}$. Pick any $a$ of the form
$$
a=\sum_{k\in\mathbb Z}r_k \;e_k\otimes e_k\qquad \mbox{ and } \quad\sum_k \mid r_k\mid^2<+\infty
$$
with all $r_k$ are different. (For instance, put $r_k=\frac{1}{\mid k\mid +1}$). Obviously, $a\in {{\sf HS}}^h$.  We affirm that there is no Hilbert-Schmidt unitary such that $SaS^*=waw^*$\end{ejem}
\begin{proof}
To prove this, suppose that there is an $w\in{\sf U} {\sf H}_{\mathbb C}$ such that $SaS^*=waw^*$. From this equation we deduce that $S^*w$ commutes with $a$, and given the particular $a$ and the fact that $S^*w$ is unitary, we have
$$
S^*w=\sum_{k\in\mathbb Z}\omega_k \;e_k\otimes e_k\qquad \mbox{ with } \mid \omega_k\mid=1
$$
because $a$ is multiplicity free. Multiplying by $S$ we get to
$$
w=\sum_{k\in\mathbb Z}\omega_k \;(Se_k)\otimes e_k=\sum_{k\in\mathbb Z}\omega_k \;e_{k+1}\otimes e_k
$$
or, in other terms, $we_k=\omega_ke_{k+1}$. Since $w$ is a compact perturbation of a scalar operator, $w$ must have a nonzero eigenvector $x$, with eigenvalue $\alpha=\mbox{e}^{i\theta}$ (since $w$ is also unitary); comparing coefficients the equation $\alpha x=wx$ reads
$$
\alpha x_k=\omega_{k-1}x_{k-1},\quad\mbox{ where }x=\sum_k x_ke_k
$$
This is impossible because $x\in l_2(\mathbb Z)$, but the previous equation leads to 
$$\mid x_k\mid=\mid x_j\mid\mbox{ for any }k,j\in\mathbb Z$$\end{proof}

\smallskip

As we see from the previous example, the two orbits do not coincide in general. For the action of the full group of unitaries we have the following:

\begin{teo}\label{full} The set $\mathfrak S_a\subset {{\sf HS}}^h$  can be given an analytic submanifold structure  if and only if the $C^*$-algebra generated by $a$ and $1$ is finite dimensional. \end{teo}
\begin{proof} The 'only if' part goes in the same lines of the proof of the previous theorem but being careful about the topologies involved, since now we must take an open set $U_a\subset {{\sf HS}}^h$ such that the map $\phi:U_a\to {\sf U} {\sf B}(H)$ is analytic. But this can be done since the polynomials $p^i_{jk}$ are now taken from $U_a$ to $\left({\sf B}(H)^n,\| \; \|_{\infty}\right)$, and the maps $+$ and $\cdot$ are analytic since $\|x.y\|_{\infty}\le  \| x\|_{_2}\|y\|_{_2}$.

The relevant part of this theorem is the 'if' part. Suppose we can prove that the orbit $\mathfrak S_a$ is closed in ${\sf B}(H)$. Then Voiculescu's theorem (see \cite{voicu}, Proposition 2.4) would tell us that $C^*(a)$ is finite dimensional. This is a deep result about $*-$representations, and the argument works in the context of ${\sf B}(H)$, but not in ${\sf H}_{\mathbb C}$ because the latter is not a $C^*$-algebra.

To prove that $\mathfrak S_a$ is closed in ${\sf B}(H)$, we first prove that it is closed in ${\sf H}_{\mathbb C}$.
To do this, observe that if $\mathfrak S_a$ is an analytic submanifold of ${{\sf HS}}^h$, then $\mathfrak S_a$ must be locally closed in the $\|\cdot\|_{_2}$ norm (in the sense that  every point $p\in \mathfrak S_a$ has an open neighborhood $U$ in ${\sf HS}^h$ such that ${\mathfrak S}_a\cap U$  is closed in $U$, see \cite{lang}). Since the action of the full unitary group is isometric, the neighborhood can be chosen uniformly, that is, there is an $\epsilon>0$ such that for all $c\in \mathfrak S_a$, the set $N_c=\{d\in \mathfrak S_a:\|c-d\|_{_2} < \epsilon\}$ is closed in the open ball $B(c,\varepsilon)=\{ x\in {\sf HS}^h: \|x-c\|< \varepsilon\}$ (with the 2-norm, of course). Now the proof that $\mathfrak S_a$ is closed in ${{\sf HS}}^h$ is straightforward, therefore we omit it. 

Now suppose $a_n=u_nau_n^*\to y$ in ${\sf B}(H)$. We claim that $\|a_n-y\|_{_2}\to 0$, which follows from a dominated convergence theorem for trace class operators (see \cite{simon}, Theorem 2.17). The theorem states that whenever $\|a_n-y\|_{\infty}\to 0$ and $\mu_k(a_n)\le\mu_k(a)$ for some $a\in {{\sf HS}}$, and all $k$ (here $\mu_k(x)$ denotes the non zero eigenvalues of $\mid x\mid$), then $\|a_n-y\|_{_2}\to 0$. 

Observe that $\mid a_n\mid =u_n \mid a\mid u_n^*$ so we have in fact equality of eigenvalues. This proves that $\mathfrak S_a$ is closed in ${\sf B}(H)$ since it is closed in ${{\sf HS}}^h$. \end{proof}

\smallskip

We proved that, when the spectrum of $a$ is finite, $S_a$ and $\mathfrak S_a$ are submanifolds of $\Sigma_1$. But more can be said: $S_a$ and $\mathfrak S_a$ are the \emph{same subset} of ${{\sf HS}}^h$ (compare with  Example \ref{noda}):
\index{operator!finite rank}

\begin{lema}\label{coincide} If $\;a\in {{\sf HS}}^h$ has finite spectrum, the orbit under both unitary groups coincide.\end{lema} 
\begin{proof}The main idea behind the proof is the fact that, when $\sigma(a)$ is finite, $a$ and $gag^*$ act on a finite dimensional subspace of $H$ (for any $g\in {\sf U}{{\sf B}(H)}$). To be more precise, let's call $S=R(a)$,  $V=R(b)$, where $b=gag^*$. Note that $V=g(S)$ so $S$ and $V$ are isomorphic, finite dimensional subspaces of $H$. Naming $T=S+V$ this is another finite dimensional subspace of $H$, and clearly $a$ and $b$ act on $T$, since they are both self-adjoint operators. For the same reason, there exist unitary operators $P,Q\in {\sf B}(T)$ and diagonal operators $D_a,D_b\in {\sf B}(T)$ such that
$$
a=PD_aP^*,\quad b=QD_bQ^*
$$
But $\sigma(b)=\sigma(gag^*)=\sigma(a)$, so $D_a=D_b:=D$. This proves that $b=QP^*aPQ^*$ (the equality should be interpreted in $T$). Now take $P_T$ the orthogonal projector in ${\sf B}(H)$ with rank $T$, and set $u=1+P_T(QP^*-1_T)P_T$. Then clearly $u\in {\sf U}{{\sf H}_{\mathbb C}}$ and $uau^*=b$.\end{proof}

\medskip

\subsection{\large Riemannian structures for the orbit $\Omega$}\label{omega}

Suppose that  there is, in fact, a submanifold structure for $\mathfrak S_a$ (resp.  $S_a$). Then the tangent map  ( $=d_1\pi_a$) has image 
$$\{va-av: v\in {\sf B}^{ah}\},
$$
where ${\sf B}$ stands for the Banach algebra ${\sf B}(H)$ (resp. ${\sf H}_{\mathbb C}$). So, in this case $$\displaystyle T_aS_a  (\mbox{ or } T_a\mathfrak S_a)=\{va-av: v\in {\sf B}^{ah}\}$$
We can go back to the manifold $\Sigma_1$ \emph{via} the usual exponential of operators; we will use the notation
$$
{\Omega}=\mbox{e}^{S_a}\quad \mbox{ or }\quad {\Omega}=\mbox{e}^{\mathfrak S_a}
$$
without further distinction. Note that $\displaystyle {\Omega}=\{u\mbox{e}^au^*\;:\; u\in {\sf U}{ {\sf B} }\}\subset \Sigma_1$ and we can identify\index{_o@$\Omega$ unitary orbit in $\Sigma_1$}
$$
T_{\mbox{e}^a}{\Omega}=\{v\mbox{e}^a-\mbox{e}^av: v\in {\sf B}^{ah}\}=\{\;i(h\mbox{e}^a-\mbox{e}^ah): h\in {\sf B}^{h}\}
$$

\begin{rem}\label{complemento} For any $p\in {\Omega}$, we have
$$
T_{p}{\Omega}=\{vp-pv: v\in {\sf B}^{ah}\}\quad\mbox{ and }\quad {T_{p}{\Omega}}^{\perp}=\{X\in {{\sf HS}}^h: [X,p]=0\}
$$
These two identifications follow from the definition of the action, and the equality
$$
\left<x,vp-pv\right>_p=4\,tr\left[(p^{-1}x-xp^{-1})v\right]
$$
\end{rem}

\begin{rem}
The submanifold $\Omega$ is connected: the curves indexed by $w\in {\sf B}^{ah}$,
$$
\gamma_w(t)=\mbox{e}^{tw}\mbox{e}^a\mbox{e}^{-tw}
$$
join $\mbox{e}^a$ to $u\mbox{e}^au^*$, assuming that $u=\mbox{e}^w$. 
\end{rem}

\smallskip

We can ask whether the curves $\gamma_w$ will be the familiar geodesics of the ambient space (equation (\ref{curvas}) of the introduction). Of course they are trivial geodesics if $a$ and $w$ commute. We will prove that this is the only case, for any $a$:

\begin{prop}\label{nosirve} For any $a\in {{\sf HS}}^h$, the curve $\gamma_w$ is a geodesic of $\;\Sigma_1$ if and only if $w$ commutes with $a$. In this case the curve reduces to the point ${\rm e}^a$.
\end{prop}
\begin{proof} The (ambient) covariant derivative for $\;\gamma_w$ (equations (\ref{covariant})  and (\ref{geo}) of the introduction) simplifies up to $w\mbox{e}^aw\mbox{e}^{-a}=\mbox{e}^aw\mbox{e}^{-a}w$ or, writing $w=ih$ ($h$ is self-adjoint)
\begin{equation}\label{esta}
\quad\quad h\mbox{e}^ah\mbox{e}^{-a}=\mbox{e}^ah\mbox{e}^{-a}h
\end{equation}\index{geodesic}
Consider the Hilbert space $(H,\left<\;,\;\right>_a)$ with inner product 
$$\left<x,y\right>_a=<{\rm e}^{-a/2}x,{\rm e}^{-a/2}y>\,$$
 where $\left<\;,\;\right>$ is the inner product of $H$. The norm of an operator $X$ is given by
$$
\|X\|_a=\sup_{\|z\|_a=1}\|Xz\|_a=\sup_{\|{\rm e}^{-a/2}z\|=1}\|{\rm e}^{-a/2}Xz\|_{\infty}=
\|{\rm e}^{-a/2}X{\rm e}^{a/2}\|_{\infty}
$$
because ${\rm e}^{-a/2}$ is an isomorphism of $H$. This equation also shows that the Banach algebras $\left({\sf B}(H),\|\cdot\|_{\infty}\right)$ and ${\sf A}=\left({\sf B}(H),\|\cdot\|_{a}\right)$ are topologically isomorphic and, as a byproduct, $\sigma_{{\sf A}}(h)\subset \mathbb R$. From the very definition it also follows easily that ${\sf A}$ is indeed a C$^*$-algebra.

A similar computation shows that $X^{*{\sf A}}={\rm e}^aX^*{\rm e}^{-a}$. Note that ${\rm e}^a$ is ${\sf A}$-self-adjoint, moreover, it is ${\sf A}$-positive. We can restate equation (\ref{esta}) as
$$
hh^{*{\sf A}}=h^{*{\sf A}}h,
$$
This equation says that $h$ is ${\sf A}$-normal, so a generalization of Weyl-von Neumann's theorem says that it can be approximated by diagonalizable operators with the same spectrum  \cite{wvn}; since $h$ has real spectrum, $h$ turns out to be ${\sf A}$-self-adjoint.  That $h$ is ${\sf A}$-self-adjoint reads, by definition, ${\rm e}^ah{\rm e}^{-a}=h^{*{\sf A}}=h$; this proves that $a$ and $h$ (and also $a$ and $w$) commute.\end{proof}

\bigskip

\subsubsection{\large $\Omega$ as a Riemannian submanifold of ${{\sf HS}}^h$}\label{plana}

We've shown earlier that the orbit of an element $a\in {{\sf HS}}^h$ has a structure of analytic submanifold of ${{\sf HS}}^h$ (which is a flat Riemannian manifold) if and only if $\Omega={\rm e}^a$ has a structure of analytic submanifold of $\Sigma_1$. 

Since the inclusion $\Omega\subset {{\sf HS}}^h$ is an analytic embedding, we can ask whether the curves
$$\gamma_w(t)=\mbox{e}^{tw}\mbox{e}^a\mbox{e}^{-tw}$$ will be geodesics of $\Omega$ as a Riemannian submanifold of ${{\sf HS}}^h$ (with the induced metric). 


We notice that the geodesic equation reads $\ddot\gamma_w(t)\perp T_{\gamma_w(t)}{\Omega}$, and we use the elementary identities $\dot\gamma=w\gamma-\gamma w$, $\ddot\gamma=w^2\gamma-2w\gamma w+\gamma w^2$; we get to the following necessary and sufficient condition using the characterization of the normal space at $\gamma(t)$ of the previous section:
$$
w^2\gamma^2-2w\gamma w\gamma +2\gamma w \gamma w- \gamma^2 w^2=0
$$
But observing that $\mbox{e}^{-wt}\gamma\;^{^{\underline{+}1}}\; \mbox{e}^{wt}=\mbox{e}^{^{\underline{+}a}}$, this equation translates into the operator condition
\begin{equation}\label{condition}
w^2\mbox{e}^{2a}-2w\mbox{e}^aw\mbox{e}^a+2\mbox{e}^aw\mbox{e}^aw-\mbox{e}^{2a}w^2=0
\end{equation}

\index{geodesic!equation}
Let's fix some notation: set $\mbox{e}^a=1+A$ with $A\in {{\sf HS}}^h$; then the tangent space at $\mbox{e}^a$ can be thought of as the subspace
$$
T_{\mbox{e}^a} {\Omega}=\{ \; i(Ah-hA):h\in {\sf B}^h\}\subset {{\sf HS}}^h
$$
and its orthogonal complement in ${{\sf HS}}^h$ is (see Remark \ref{complemento})
$$
{T_{\mbox{e}^a} {\Omega}}^{\perp}=\{ \; X\in {\sf B}^h: [X,A]=0\}
$$
It should be noted that both subspaces are closed by hypothesis. Then equation (\ref{condition}) can be restated as
\begin{equation}\label{coso}
h^2A^2-2hAhA+2AhAh-A^2h^2=0
\end{equation}
where $h$ is the hermitian generating the curve 
$$\gamma(t)=1+\mbox{e}^{ith}A\mbox{e}^{-ith}=\mbox{e}^{ith}\mbox{e}^a\mbox{e}^{-ith}$$

\bigskip

Let's consider the case when $A^2=A$:

\index{operator!projector}
\begin{rem}\label{igual} If $A^2=A$, $A$ must be a finite rank orthogonal projector (since $A={\rm e}^a-1$ and $ a$ is a Hilbert-Schmidt operator).  Hence $\sigma(a)$ consists of two points, and in this case (Remark \ref{coincide}) the orbit with the full unitary group and the orbit with the Hilbert-Schmidt  unitary group are the same set.
\end{rem}

\smallskip

Observe that when $A$ is a projector, we have a matrix decomposition of the tangent space of $\Sigma_1$, namely ${{\sf HS}}^h=A_0\oplus A_1$, where 
$$
A_0=\left\{ \left(  \begin{array}{cc} x_{11}  & 0   \\  0 & x_{22}  \\  \end{array}   \right)  \right\}\quad\mbox{ and }\quad A_1=\left\{ \left(  \begin{array}{cc} 0 & x_{12}    \\   x_{21} & 0  \\  \end{array}   \right)  \right\}
$$

In this decomposition, $x_{11}=AhA\;,\; x_{22}=(1-A)h(1-A)$ are self-adjoint operators (since $h$ is) and also $x_{12}^*=x_{21}=(1-A)hA$ for the same reason.

\begin{teo}\label{cartan} Whenever $A={\rm e}^a-1$ is a projector, any curve of the form $ \gamma(t)=\mbox{e}^{ith}\mbox{e}^a\mbox{e}^{-ith}$ with $h$ self-adjoint and co diagonal is a geodesic of $\;\Omega\subset {{\sf HS}}^h$ 
\end{teo}
\begin{proof} 

Note that $A_0={T_{\mbox{e}^a} {\Omega}}^{\perp}$, and $A_1={T_{\mbox{e}^a} {\Omega}}$; 
note also that equation (\ref{coso}) translates in this context to $x_{11}x_{12}=x_{12}x_{22}$, a condition which is obviously fulfilled by $h\in A_1$.\end{proof}

\begin{rem}\label{remarc} Equation (\ref{coso}) translates exactly in '$h_0$ commutes with $h_1$' whenever $h=h_0+h_1\in {{\sf HS}}^h$, and we have
$$
[A_0,A_1]\subset A_1\quad [A_0,A_0]\subset A_0\quad [A_1,A_1]\subset A_0
$$
Since the orbit under both unitary groups coincide (Remark \ref{igual}), assume that we are acting with $G={\sf U}{ \sf B }$; since the tangent space at the identity of this group can be identified with ${\sf B}^{ah}$, the above commutator relationships say that $iA_0\oplus iA_1$ is a Cartan decomposition of the Lie algebra $\mathfrak g={\sf B}^{ah}$. It is apparent that $iA_0$ is the vertical space, and $iA_1$ is the horizontal space. Moreover, \index{Cartan decomposition}\index{graded algebra}\index{vertical space}
$$
A_0\cdot A_0\subset A_0\qquad A_1\cdot A_1\subset A_0\qquad A_0\cdot A_1\subset A_1\qquad A_1\cdot A_0\subset A_1
$$
\end{rem}

\begin{coro}\label{no} If ${\rm e}^a-1$ is an orthogonal projector, there is no point $p\in\Omega$ such that $\Omega$ is geodesic at $p$.
\end{coro}

\begin{rem}
In the paper \cite{cpr93} by Corach, Porta and Recht, the authors study the differential geometry of self-adjoint projections in a $C^*$-algebra. The authors show the role of the graded decomposition of the algebra in the characterization of the geodesics for the Finsler structure that this space carries; the geodesics we obtained are similar to the ones obtained in that paper. 
\end{rem}

\subsubsection{\large $\Omega$ as a Riemannian submanifold of $\Sigma_1$}\label{posta}

In this section we give $\Omega$ the induced Riemannian metric as a submanifold of $\Sigma_1$, and discuss shortly the induced exponential map.

\smallskip

Recall that covariant derivative in the ambient space is given by $\nabla_{\dot\gamma}\dot\gamma=\ddot\gamma-\dot\gamma\gamma^{-1}\dot\gamma\;$ and the orthogonal space to $p\in \Omega$ are the operators commuting with $p$, so $\nabla_{\dot\gamma}\dot\gamma\perp T_{\gamma}\Omega$ if and only if
\begin{equation}\label{sub}
\ddot\gamma\gamma-\gamma\ddot\gamma+\gamma\dot\gamma\gamma^{-1}\dot\gamma-\dot\gamma\gamma^{-1}\dot\gamma\gamma=0
\end{equation}
This is an odd equation; we know that any curve in $\Omega$ starting at $p={\rm e}^a$ must be of the form $\gamma(t)=g(t){\rm e}^ag(t)^*$ for some curve of unitary operators $g$.
\index{geodesic!equation}

\smallskip

For the particular curves $\gamma(t)={\rm e}^{ith}{\rm e}^a{\rm e}^{-ith}$, $h(t)=ith$, so  $\dot h(t)=ih$, and $\ddot h(t)\equiv 0$; equation (\ref{sub}) reduces to the operator equation
\begin{equation}\label{larga}
h{\rm e}^ah{\rm e}^{-a}+ h {\rm e}^{-a}h {\rm e}^a=  {\rm e}^{-a}h{\rm e}^ah+{\rm e}^{a} h{\rm e}^{-a}h
\end{equation}
or $X^*=X$, where $X=h{\rm e}^{a}h {\rm e}^{-a}+h{\rm e}^{-a}h{\rm e}^a$.

\medskip

Recall that the unitary groups ${\sf U}{{\sf B}(H)}$ and ${\sf U}{{\sf H}_{\mathbb C}}$ induce the same manifold $\Omega\subset\Sigma_1$ when the spectrum of $a$ is finite.  Throughout $[\; ,\; ]$ stands for the usual commutator of operators.

\begin{teo}\label{ultimo}
Assume ${\rm e}^a=1+A$ with $A$ an orthogonal projector, and $\Omega\subset\Sigma_1$ is the unitary orbit of ${\rm e}^a$. Then 
\begin{enumerate}
\item[(1)] $\Omega$ is a Riemannian submanifold of $\;\Sigma_1$.
\item[(2)] $T_p\Omega=\{  i[x,p]:x\in {{\sf HS}}^h\}$ and $T_p\Omega^{\perp}=\{  x\in {{\sf HS}}^h: [x,p]=0\}$.
\item[(3)] The action of the unitary group is isometric, namely
$$\mbox{\rm dist}^{\Omega}\left(  upu^*,uqu^*   \right)=\displaystyle\mbox{\rm dist}^{\Omega}\left(  p,q    \right)$$ for any unitary operator $u\in {\sf B}(H)$.
\item[(4)] For any $v=i[x,p]\in T_p\Omega$, the exponential map is given by  
$$\;{\rm exp}_{p}^{\Omega}(v)={\rm e}^{ighg^*}p\,{\rm e}^{-ighg^*}$$ where $p=g{\rm e}^ag^*$ and $h$ is the co diagonal part  of $g^*xg$ (in the matrix representation of Proposition \ref{cartan}). In particular, the exponential map is defined in the whole tangent space.
\item[(5)] If $p=g{\rm e}^ag^*$, $q=w{\rm e}^aw^*$, and $h$ is a self-adjoint, co diagonal operator such that $w^*g{\rm e}^{ih}$ commutes with $e^a$, then the curve $\gamma(t)={\rm e}^{itghg^*}p{\rm e}^{-itghg^*}$ is a geodesic of $\;\Omega\subset\Sigma_1$, which joins $p$ to $q$.
\item[(6)] If we assume that $h\in {{\sf HS}}^h$, then $L(\gamma)=\frac{\sqrt{2}}{2}\; \|h\|_{_2}$
\item[(7)] The exponential map $\;{\rm exp}_{p}^{\Omega}:T_p\Omega\to \Omega$ is surjective.\end{enumerate}
\end{teo}
\begin{proof} 
\index{group!of isometries}
Statements (1) and (2) are a consequence of Remark \ref{igual} and Theorems \ref{hs} and \ref{full}. 

Statement (3) is obvious because the action of the unitary group is isometric for the 2-norm (see Lemma \ref{invar}). 

To prove statement (4), take $x\in {\sf HS}^h$, and set 
$$v=i[x,p]=i(xgAg^*-gAg^*x)=ig[g^*xg,{\rm e}^a]g^*$$
Observe that
$${\rm e}^{-a}=(1+A)^{-1}=1-\frac{1}{2}A$$
Rewriting equation (\ref{larga}), we obtain 
$$h^2A-Ah^2+2AhAh-2hAhA=0$$
Now if $y=g^*xg$, take $h=$ the co diagonal part of $y$; clearly $hA-Ah=yA-Ay$, so 
$$\gamma_1(t)={\rm e}^{ith}{\rm e}^a{\rm e}^{-ith}$$ is a geodesic of $\Omega$ starting at $r={\rm e}^a$ with initial speed $w=i[y,{\rm e}^{a}]=g^*vg$ (see Proposition \ref{cartan}).
Now consider $\gamma=g\gamma_1g^*$. Clearly $\gamma$ is a geodesic of $\Omega$ starting at $p=g{\rm e}^ag^*$ with initial speed $v$.

To prove (5), note that
$$\gamma(t)=g{\rm e}^{iht}{\rm e}^a{\rm e}^{-iht}g^*={\rm e}^{itghg^*}g{\rm e}^ag^*{\rm e}^{itghg^*}={\rm e}^{itghg^*}p{\rm e}^{itghg^*}$$
which shows that $\gamma(0)=p$ and $\gamma(1)=q$ because $w^*g{\rm e}^{ih}{\rm e}^{a}={\rm e}^{a}w^*g{\rm e}^{ih}$.

To prove (6), we can assume that $p={\rm e}^a$, and then
$$
L(\gamma)^2=\|[h,p]\|_p^2=\|[h,{\rm e}^a]\|_{{\rm e}^a}^2=4\cdot tr(2h{\rm e}^ah{\rm e}^{-a}-2h^2)
$$
Now write $h$ as a matrix operator $[0,Y^*,Y,0]\in A_1$ (see Proposition \ref{cartan}), to obtain 
$$tr(2h{\rm e}^ah{\rm e}^{-a}-2h^2)=tr(Y^*Y)=\textstyle\frac12\, tr(h^2),$$
 hence $L(\gamma)^2=2\, tr(h^2)= \frac12 \|h\|_{_2}^{^2}$ as stated. 

The assertion in (7) can be deduced from folk results (see \cite{brown}) because $q=w{\rm e}^a w^*$ and $p=g {\rm e}^ag^*$ are finite rank projectors acting on a finite dimensional space (see the proof 
of Lemma \ref{coincide}). \end{proof}

\bigskip

\section{Concluding remarks}

\begin{rem}
Theorem \ref{hs} does not answer whether is it necessary that the spectrum of $a$ should be finite for the orbit to be a submanifold, when we act with ${\sf U}({{\sf H}_{\mathbb C}})$ (see Remark \ref{nose}). The problem can be stated in a more general form:
\begin{itemize}
\item Choose any involutive Banach algebra with identity ${\sf B}$, take $a\in {\sf B}$ such that $a^*=a$, and denote ${\sf U}_{ {\sf B} }=\{u\in {\sf B}:u^*=u^{-1}\}$, the unitary group of ${\sf B}$.
\item Name $S_a$ the image of the map $\pi_a:{\sf U}_{ {\sf B} }\to{\sf B}$ which assigns $u\mapsto uau^*$
\item Is the condition "$a$ has finite spectrum" necessary for the set $S_a\subset {\sf B}$ to be closed?
\end{itemize}
\end{rem}

\bigskip

\noindent
Gabriel Larotonda\\
Instituto de Ciencias \\
Universidad Nacional de Gral. Sarmiento \\
J. M. Gutierrez 1150 \\
(1613) Los Polvorines \\
Argentina  \\
e-mail: glaroton@ungs.edu.ar

\end{document}